\documentclass{amsproc}

\usepackage{times}
\usepackage{amssymb}
\usepackage{amsmath}
\usepackage{aspmproc}

\begin{document}

\newtheorem{theorem}{Theorem}[section]
\newtheorem{corollary}[theorem]{Corollary}
\newtheorem{lemma}[theorem]{Lemma}
\newtheorem{proposition}[theorem]{Proposition}
\newtheorem{definition}[theorem]{Definition}
\newtheorem{remark}{Remark}

\newcommand{\A}{{\bf \mathcal A}}
\newcommand{\Ab}{{\mathbb A}}
\newcommand{\Alift}{{A^{\wedge}}}
\newcommand{\Aoplift}{{{\mathcal A}^{\wedge}}}
\newcommand{\B}{{ \bf \mathcal B }}
\newcommand{\Bb}{{ \mathbb B }}
\newcommand{\C}{{\mathbb C}}
\newcommand{\D}{{\bf \mathcal D}}
\newcommand{\E}{{ \Bbb E }}
\newcommand{\F}{{ \mathcal F }}
\newcommand{\1}{{\bf 1}}

\newcommand{\HH}{{\mathcal H}}
\newcommand{\PP}{{\mathbb P}}
\newcommand{\Q}{{\mathbb Q}}
\newcommand{\LL}{{\mathcal L}}
\newcommand{\R}{{\Bbb R}}
\newcommand{\X}{{\mathbb X}}
\newcommand{\Y}{{\mathbb Y}}
\newcommand{\Z}{{\mathbb Z}}

\newcommand{\trace}{{\rm trace}}
\def\paral{/\kern-0.3em/}
\def\parals_#1{/\kern-0.3em/_{\!#1}}
\renewcommand{\thefootnote}{}

\title{Equivariant Diffusions on Principal Bundles }
\author{K. D. Elworthy \and Yves Le Jan \and Xue-Mei Li}
\maketitle

\noindent
\footnote{Research supported by  EPSRC grants GR/NOO 845 and GR/S40404/01,
  NSF grant DMS 0072387, and EU grant ERBF MRX CT 960075 A. }

 Let $\pi: P\to M$ be a smooth principal bundle with
structure group $G$. This means that there is a $C^\infty$ right
multiplication $P\times G\to P$, $u\mapsto u\cdot g$ say, of the
Lie group $G$ such that $\pi$ identifies the space of orbits of
$G$ with the manifold $M$ and $\pi$ is  locally trivial in the
sense that each point of $M$ has an open neighbourhood $U$ with a
diffeomorphism

\begin{picture}(190,60)(-20,18)
\put(25,58){$\tau_U: \pi^{-1}(U)$}
 \put(87, 60){\vector(1,0){52}}
 \put(153, 58){$U\times G$}
\put(60,50){\vector(2,-1){37}} \put(106, 22){$U$}
\put(161,50){\vector(-2,-1){37}}
\end{picture}\\
over $U$, which is equivariant with respect to the right action of
$G$, i.e. if $\tau_u(b)=(\pi(b),k)$ then $\tau_u(b\cdot
g)=(\pi(b),kg)$. Assume for simplicity that $M$ is compact. Set
$n=dim M$. The fibres, $\pi^{-1}(x)$, $x\in M$ are diffeomorphic
to $G$ and their tangent spaces $VT_uP(=ker T_u \pi)$, $u\in P$,
are the `vertical' tangent spaces to $P$. A {\it connection} on
$P$, (or on $\pi$) assigns a complementary `horizontal' subspace
$HT_uP$ to $VT_uP$ in $T_uP$ for each $u$, giving a smooth
horizontal sub-bundle $HTP$ of the tangent bundle $TP$ to $P$.
Given such a connection it is a classical result that for any
$C^1$ curve: $\sigma: [0, T]\to M$ and $u_0\in
\pi^{-1}(\sigma(0))$ there is a unique horizontal
 $\tilde \sigma: [0, T]\to P$ which is a lift of $\sigma$, i.e.
 $\pi (\tilde \sigma(t))=\sigma(t)$ and has
$\tilde \sigma(0)=u_0$.

In his startling ICM article \cite{Ito} It\^o showed how this construction
could be extended to give horizontal lifts of the sample paths of diffusion
processes. In fact he was particularly concerned with  the case when $M$ is
 given a Riemannian metric $\langle, \rangle_x$, $x\in M$, the diffusion
is Brownian motion on $M$, and $P$ is the orthonormal frame bundle
$\pi: OM\to M$. Recall that each $u\in OM$ with $u\in \pi^{-1}(x)$ can be
considered as an isometry $u: \R^n \to T_xM$,  $\langle, \rangle_x$ and a
 horizontal lift $\tilde \sigma$ determines parallel translation of tangent
 vectors along $\sigma$
\begin{eqnarray*}
\parals_t\equiv \paral(\sigma)_t:&&  T_{\sigma(\cdot)}M\to T_{\sigma(t)}M\\
&&v\mapsto \tilde \sigma(t)(\tilde \sigma(0))^{-1}v.
\end{eqnarray*}
The resulting parallel translation along Brownian paths extends also to
 parallel translation of forms and elements of $\wedge^p TM$. This enabled
It\^o to use his construction to obtain a semi-group acting on
differential forms
$$P_t\phi={\mathbb E} (\parals_t^{-1})_*(\phi)={\mathbb E} \phi(\parals_t-).$$
 As he pointed out this is not the semi-group
generated by the Hodge-Kodaira Laplacian, $\Delta$. To obtain that
generated by the Hodge-Kodaira Laplacian, $\Delta$, some modification
had to be made since the latter contains
zero order terms, the so called Weitzenbock curvature terms.
The resulting probabilistic expression for the heat semi-groups
on forms has played a major role in subsequent development.

In \cite{Elworthy-LeJan-Li-principal}
 we go in the opposite direction starting  with a diffusion with
 smooth generator $\mathcal B$ on
$P$, which is $G$-invariant and so projects to a diffusion
generator $\mathcal A$ on $M$. We assume the symbol
$\sigma_{\mathcal A}$ has constant rank so determining a
sub-bundle $E$ of $TM$, (so $E=TM$ if $\A$ is elliptic). We show
that this set-up induces a `semi-connection' on $P$ over $E$ (a
connection if $E=TM$) with respect to which $\B$ can be decomposed
into a horizontal component $\A^H$ and a vertical part $\B^V$.
Moreover any vertical diffusion operator such as $\B^V$ induces
only zero order operators on sections of associated vector
bundles.

There are two particularly interesting  examples. The first when
$\pi: GLM\to M$ is the full linear frame bundle and we are given a stochastic
 flow $\{\xi_t: 0\le t<\infty\}$ on $M$, generator $\A$,
inducing the diffusion $\{u_t: 0\le t<\infty\}$ on GLM by
$$u_t=T\xi_t(u_0).$$
Here we can determine the connection on GLM in terms of the
LeJan-Watanabe connection of the flow  \cite{LeJan-Watanabe},
\cite{Accardi-Mohari}, as defined in
\cite{Elworthy-LeJan-Li-Tani}, \cite{Elworthy-LeJan-Li-book},
 in particular giving conditions when it is a Levi-Civita connection.
 The zero order operators arising  from the vertical components, can be
 identified with generalized Weitzenbock curvature terms.

The second example slightly extends the above framework by letting
 $\pi: P\to M$ be the evaluation map
 on the diffeomorphism  group Diff$M$ of $M$ given by
 $\pi(h):=h(x_0)$ for a fixed point  $x_0$ in $M$.
 The  group $G$ corresponds to the group of diffeomorphisms
 fixing $x_0$. Again we take a flow $\{\xi_t(x): x\in M,  t\ge 0\}$
on $M$, but now the process on Diff$M$ is just the right invariant
  process determined by $\{\xi_t: 0\le t<\infty\}$.
In this case the horizontal lift to the diffeomorphism group
 of the diffusion
 $\{\xi_t(x_0): 0\le t<\infty\}$ on $M$ is obtained by `removal of
 redundant noise', c.f. \cite{Elworthy-LeJan-Li-book} while the
vertical process is a flow of diffeomorphisms preserving $x_0$, driven by
 the redundant noise.

Here we report briefly on some of the main results to appear in
\cite{Elworthy-LeJan-Li-principal} and give details of a more
probabilistic version Theorem \ref{theorem-kernel-decomposition}
below: a skew product decomposition which, although it has a
statement not explicitly mentioning connections, relates to
It\^o's pioneering work on the existence of horizontal lifts. The
derivative flow example and a simplified version of the stochastic
flow example are described in \S~\ref{section-stochastic-flows}.

The decomposition and lifting apply in much more generality than
with the full structure of a principal bundle, for example to
certain skew products and invariant processes on foliated
manifolds. This will be reported on later. Earlier work on such
decompositions includes
 \cite{Elworthy-Kendall}  \cite{Liao89}.

\section{Construction}\label{Construction}

{\bf A. }
If ${\mathcal A}$ is a second order differential operator on a manifold $X$,
 denote by $\sigma^{\mathcal{A}}: T^*X\to TX $ its symbol determined by
\[
df\left( \sigma^{\mathcal{A}}(dg)\right) =\frac 12{\mathcal A}\left(
fg\right) -\frac 12{\mathcal A}\left( f\right) g-\frac 12f{\mathcal A}\left(
g\right),
\]
for $C^2$ functions  $f, g$.
The operator is said to be {\it semi-elliptic}  if
 $\displaystyle{df\left( \sigma ^{\mathcal A}(df)\right) \ge 0}$
 for each $\displaystyle{f\in C^2(X)}$, and {\it elliptic} if the
inequality holds strictly. Ellipticity is equivalent to
 $\displaystyle{\sigma^{\mathcal{A}}}$ being  onto. It is called a
{\it diffusion operator} if it is semi-elliptic and annihilates
constants, and is {\it smooth} if it sends smooth functions to
smooth functions.

\medskip

Consider a smooth map $p:N\to M$ between smooth manifolds $M$ and $N$.
By a {\it lift}  of a diffusion operator ${\mathcal  A}$ on $M$ over
  $p$ we mean a diffusion operator ${\mathcal  B}$ on $N$ such that
\begin{equation}\label{operator-lift}
{\mathcal  B}(f\circ p)=({\mathcal  A}f)\circ p
\end{equation}
for all $C^2$ functions $f$ on $M$. Suppose ${\mathcal  A}$ is a
smooth diffusion operator on $M$ and ${\mathcal  B}$ is a lift of
${\mathcal  A}$.

\begin{lemma}
\label{le:commutative} \index{commutative} Let $\sigma
^{\mathcal{B}}$ and $\sigma ^{\mathcal{A}}$ be respectively the
symbols for $\mathcal{B}$ and $\mathcal{A}$. The following diagram
is commutative for all $u\in p ^{-1}(x)$, $x\in M$:

\begin{picture}(200, 80)(-40,-57)
\put(30,0) {$T_u ^*N$} \put(85,8){$\sigma_u^{\mathcal B}$}
\put(55,2){\vector(1,0){70}} \put(137,0){$T_uN$}
\put(45,-45){\vector(0,1){40}} \put(18,-25){$(Tp)^*$}
\put(32,-55){$T_x ^*M$} \put(135,-4){\vector(0,-1){40}}
\put(133,-55){$T_xM$.} \put(65,-52){\vector(1,0){60}} \put(87,
-62){$\sigma_x^{\mathcal A}$} \put(141,-25){$Tp$}
\end{picture}
\end{lemma}

{\bf B. Semi-connections on principal bundles. }
Let $M$ be a smooth finite dimensional manifold and
 $\displaystyle{P\left( M,G\right)} $ a
 principal fibre bundle over $M$ with structure group $G$ a  Lie group.
 Denote by  $\displaystyle{\pi :P\rightarrow M}$  the projection and $R_a$
the right  translation by $a$.

\begin{definition}
Let $E$ be a sub-bundle of $TM$ and $\pi:P\to M$ a principal $G$-bundle.
An $E$ semi-connection on $\pi:P\to M$ is a smooth sub-bundle
$H^ETP$ of $TP$ such that

\begin{enumerate}
\item[(i)]
$T_u\pi$ maps the fibres $H^ET_uP$ bijectively onto $E_{\pi(u)}$
for all $u\in P$.
\item[(ii)]
$H^ETP$ is $G$-invariant.
\end{enumerate}
\end{definition}

\noindent{\bf Notes. }\\
 (1) Such a semi-connection determines and
is determined by, a smooth horizontal lift:
$$h_u: E_{\pi(u)}\to T_uP, \qquad u\in P$$ such that
\begin{enumerate}
\item
[(i)] $T_u\pi \circ h_u(v)=v$, for all $v\in E_x\subset T_xM$;
\item[(ii)] $h_{u\cdot a}=T_uR_a\circ h_u$.
\end{enumerate}
The {\it horizontal
subspace} $H^ET_uP$ at $u$ is then the image at $u$ of $h_u$, and
the composition $h_u\circ T_uP$ is a projection onto $H^ET_uP$.\\
(2)
 Let $F=P\times V/\sim$ be an associated vector bundle to $P$ with
 fibre $V$. An element of $F$ is an equivalence class $[(u,e)]$ such
that $(ug, g^{-1}e)\sim (u,e)$. Set $\tilde u(e)=[(u,e)]$.
 An $E$ semi-connection on $P$ gives a covariant derivative
on $F$. Let $Z$ be a section of $F$ and $w\in E_x\subset T_xM$,
the covariant derivative $\nabla_w Z\in F_x$  is defined, as usual
for
 connections, by
$$\nabla_w Z=u(d\tilde Z(h_u(w)), \qquad u\in \pi^{-1}(x)=F_x. $$
 Here $\tilde Z: P\to V$ is
 $\tilde Z(u)={\tilde u}^{-1}Z\left(\pi(u)\right)$ considering $\tilde u$
as an isomorphism $\tilde u:V\to F_{\pi(u)}$. This agrees with the
`semi-connections on $E$' defined in Elworthy-LeJan-Li
\cite{Elworthy-LeJan-Li-book}
 when $P$ is taken to be the linear frame bundle of $TM$  and $F=TM$.
As described there, any semi-connection can be completed to a genuine
 connection, but not canonically.
\medskip

 Consider on $P$ a diffusion generator  ${\mathcal B}$, which is {\it equivariant}, i.e.
$${\mathcal B}f\circ R_a={\mathcal B}(f\circ R_a), \qquad
\forall f,g \in C^2(P,R),\; a\in G.$$
 The operator ${\mathcal B}$ induces an
operator ${\mathcal A}$ on the base manifold $M$ by setting
 \begin{equation}
\label{operator-equivariant-induced}
{\mathcal A}f(x)={\mathcal B}\left( f\circ \pi \right) (u),
\hskip 20pt u\in \pi^{-1}(x), f\in C^2(M),
\end{equation}
which  is well defined since
 $${\mathcal B}\left( f\circ \pi\right) \left( u\cdot a\right)
={\mathcal B}\left( \left( f\circ \pi \right) \right) (u).$$

Let $E_x:=Image(\sigma_x^{\mathcal  A})\subset T_xM$, the image of
 $\sigma_x^{\mathcal  A}$. Assume the dimension of $E_x=p$, independent
 of $x$. Set $E=\cup_x E_x$. Then $\pi: E\to M$ is a sub-bundle of $TM$.

\begin{theorem}
\label{theorem:connection}
Assume $\sigma^{\mathcal A}$ has constant rank. Then $\sigma^{\mathcal B}$
 gives rise to a semi-connection on  the principal bundle $P$ whose
 horizontal map is given by
 \begin{equation}
\label{horizontal-lift}
h_u(v)=\sigma^{\mathcal B}\left( (T_u\pi )^* \alpha\right)
\end{equation}
where $\alpha\in T_{\pi (u)}^*M$  satisfies $\sigma^{\mathcal
A}_x(\alpha)=v$.
\end{theorem}

\noindent
{\it Proof. }
 To prove $h_u$ is well defined we only need to show
$\psi(\sigma^{\mathcal B} (T_u \pi ^*(\alpha)))=0$ for every 1-form
$\psi$ on $P$ and for every  $\alpha$ in ker $\sigma_x^{\mathcal
A}$. Now  $\sigma^{\mathcal A}\alpha=0$  implies  by
Lemma \ref{le:commutative} that
$$0=\alpha \sigma^{\mathcal A} (\alpha)
=(T\pi) ^*(\alpha) \sigma^{\mathcal B} ((T\pi )^*(\alpha)).$$

Thus  $T\pi^*(\alpha)\sigma^{\mathcal B} (T\pi ^*(\alpha))=0$.
On the other hand we may consider
 $\sigma^{\mathcal B}$ as a bilinear  form on $T^*P$ and then for all
$\beta\in T_u^*P$,
\begin{eqnarray*}
&&\sigma^{\mathcal B} (\beta+t (T\pi)^*(\alpha), \beta+t(T\pi)^*(\alpha))\\
&&=\sigma^{\mathcal B}(\beta, \beta) +
 2t \sigma^{\mathcal B}(\beta,  (T\pi)^*(\alpha))+
t^2\sigma^{\mathcal B} ((T\pi)^*\alpha, (T\pi)^* \alpha)\\
&&=\sigma^{\mathcal B}(\beta, \beta)
+2t \sigma^{\mathcal B}(\beta, (T\pi)^*(\alpha)).
\end{eqnarray*}
Suppose
 $\displaystyle{\sigma^{\mathcal B} (\beta, (T\pi)^*(\alpha))\not = 0}$.
 We can then choose $t$ such that
$$\displaystyle{\sigma^{\mathcal B}
(\beta+t(T\pi)^*(\alpha), \beta+t (T\pi)^*(\alpha))<0},$$
 which contradicts  the  semi-ellipticity of  ${\mathcal B}$.

We must verify (i) $T_u\pi \circ h_u(v)=v$,
 $v\in E_x\subset T_xM$ and (ii) $h_{u\cdot a}=T_uR_a\circ h_u$.
The first is immediate by Lemma \ref{le:commutative} and for the
second use the fact that   $T\pi \circ TR_a=T\pi $ for all $a\in
G$ and the equivariance of $\sigma^{\mathcal B}$.
\hfill\rule{3mm}{3mm}

\section[Horizontal lift of diffusions]
{Horizontal lifts of diffusion operators and decompositions of
 equivariant operators}

{\bf A.}
 Denote by $C^\infty \Omega^p$ the  space of smooth  differential p-forms
  on a manifold $M$. To each diffusion operator ${\mathcal A}$ we shall
associate a unique operator $\delta^{\mathcal A}$. The horizontal lift of
${\mathcal A}$ can be defined to be the unique operator such that the
associated   operator  $\bar \delta$ vanishes on vertical  1-forms and
such that   $\bar \delta$ and $\delta^\A$ are intertwined by the lift map
$\pi^*$ acting on 1-forms.

\begin{proposition}
\label{proposition-delta} For each smooth diffusion operator
$\mathcal{A}$ there is a unique  smooth differential operator
$\displaystyle{\delta ^{\mathcal{A}}: C^\infty(\Omega^1) \to
C^\infty \Omega^0}$ such that
\begin{enumerate}
\item[(1)]  $\delta ^{\mathcal{A}}\left( f\phi \right) =df\sigma ^{\mathcal{A}}( \phi ) _x
+f \cdot \delta ^{\mathcal{A}}\left( \phi \right) $
\item[(2)]  $\delta ^{\mathcal{A}}\left( df\right) ={\mathcal A}(f) .$
\end{enumerate}
\end{proposition}

For example if ${\mathcal A}$ has H\"ormander representation
\[
{\mathcal A}={1\over 2}\sum_{j=1}^m{\mathcal L}_{X^j}{\mathcal
L}_{X^j} + {\mathcal L}_A \] for some $C^1$ vector fields $X^i$,
$A$ then
$$\delta ^{\mathcal{A}}
  = {1\over 2} \sum_{j=1}^m {\mathcal L}_{X^j}\iota_{X^j}+\iota _A$$
where $\iota_{A}$ denotes the interior product of the vector field
 $A$ acting on differential forms.

\begin{definition}
Let $S$ be a $C^\infty$ sub-bundle of $TN$ for some smooth manifold $N$. A
diffusion operator ${\mathcal  B}$ on $N$ is said to be {\it along $S$}
if $\delta^{\mathcal  B} \phi=0$ for all 1-forms $\phi$ which vanish on $S$;
it is said to be strongly cohesive if $\sigma^{\B}$
has constant rank and $\B$ is along the image of $\sigma^{\B}$.
\end{definition}

To be along $S$ implies  that any H\"{o}rmander form representation of $\B$ uses only vector fields 
which are sections of $S$.

\medskip

\begin{definition}
When a diffusion operator ${\mathcal  B}$ on $P$ is along
  the vertical foliation VTP of the $\pi: P\to M$
  we say ${\mathcal  B}$ is {\it vertical},  and when
 the bundle has a semi-connection and ${\mathcal  B}$ is along the horizontal
 distribution we say ${\mathcal  B}$ is horizontal.
\end{definition}

If $\pi: P\to M$ has an $E$ semi-connection and $\A$ is a smooth
diffusion operator along $E$ it is easy to see that $\A$ has a
unique horizontal lift $\A^H$, i.e. a smooth diffusion operator
$\A^H$ on $P$ which is horizontal and is a lift of $\A$ in the
sense of (\ref{operator-lift}). By uniqueness it is equivariant.

\bigskip
 {\bf B.} The action of $G$ on $P$ induces a
homomorphism of the Lie algebra
 ${\mathfrak g}$ of $G$ with the algebra of right invariant vector fields
on $P$: if $\alpha\in {\mathfrak g}$,
$$A^\alpha(u)=\left.{d\over dt}\right|_{t=0}\; u \exp(t\alpha),$$
and $A^\alpha$ is called the fundamental vector field corresponding to
$\alpha$. Take  a basis  ${A_1, \dots, A_k}$ of ${\mathfrak
g}$ and denote the corresponding fundamental vector fields by
 $\{A_i^*\}$.

We can now give one of the main results from
\cite{Elworthy-LeJan-Li-principal}:
 \begin{theorem}
\label{theorem-operator-decomposition}
 Let ${\mathcal B}$ be an equivariant  operator on $P$ with
 ${\mathcal A}$  the induced operator on the base manifold.
 Assume $\A$ is strongly cohesive. Then
 there is a unique
 semi-connection on $P$ over $E$ for which $\B$ has a decomposition
$$\B=\A^H+\B^V,$$
where $\A^H$ is horizontal and $\B^V$ is vertical. Furthermore
${\mathcal B}^V$ has the  expression
 $\sum \alpha^{ij} {\mathcal  L}_{A_i^*}{\mathcal  L}_{A_j^*}
  +\sum\beta^k{\mathcal  L}_{A_k^*}$,
 where $\alpha^{ij}$ and $\beta^k$ are smooth functions on $P$, given by
$\displaystyle{\alpha^{k\ell}
=\tilde \omega^k\left(\sigma^{\mathcal  B}(\tilde \omega ^\ell)\right)}$,
and $\displaystyle{\beta^\ell=\delta^{\mathcal  B}(\tilde \omega^\ell)}$
for  $\tilde\omega$ any connection 1-form on $P$ which vanishes on
the horizontal subspaces of this semi-connection.
 \end{theorem}
We shall only prove the first part of Theorem
\ref{theorem-operator-decomposition} here.
 The semi-connection is the one given by Theorem \ref{theorem:connection},
and we define $\A^H$ to be the horizontal lift of $\A$. The proof that
$\B^V:=\B-\A^H$ is vertical
 is simplified by using the fact that a diffusion operator
$\D$ on $P$ is vertical if and only if for all $C^2$ functions
$f_1$ on $P$ and $f_2$ on $M$
 \begin{equation}
\label{verticality}
\D(f_1(f_2\circ \pi))=(f_2\circ \pi)\D(f_1).\end{equation}

Set $\tilde f_2=f_2\circ \pi$.  Note
 $$\left({\mathcal B}-{\mathcal A}^H \right)(f_1\tilde f_2)
 =  \tilde f_2({\mathcal B}-{\mathcal A}^H)f_1
    +f_1 ({\mathcal B}-{\mathcal A}^H)\tilde f_2
  +2 (df_1)\sigma^{ {\mathcal B}-{\mathcal A}^H}(d\tilde f_2).$$
Therefore to show
 $({\mathcal B}-{\mathcal A}^H)$ is vertical we only need to prove
 $$\displaystyle{f_1 ({\mathcal B}-{\mathcal A}^H)\tilde f_2
 +2(df_1)\sigma^{ {\mathcal B}-{\mathcal A}^H}(d\tilde f_2)=0}.$$
Recall Lemma \ref{le:commutative}  and use the natural extension
of $\sigma^{\A}$ to $\sigma^{\A}: E^*\to E$ and the fact that
by (\ref{horizontal-lift}) $h\circ \sigma^\A_x=\sigma^{\B}(T_u \pi)^*$
to see
 \begin{eqnarray*}
 \sigma^{{\mathcal A}^H} (d \tilde f_2)
&=&\left(h\circ \sigma^{\mathcal A}h^* \right)(df_2 \circ T\pi)
=h \circ \sigma^{\mathcal A} df_2 \\
 &=&\sigma^{\mathcal  B}(df_2\circ T\pi )
=\sigma^{\mathcal  B}(d\tilde f_2),
 \end{eqnarray*}
and so $\sigma^{({\mathcal B}-{\mathcal A}^H)}(d \tilde f_2)=0$.
Also by equation (\ref{operator-lift})
$$({\mathcal B}-{\mathcal A}^H)\tilde f_2
={\mathcal A}f_2 -{\mathcal A}^H \tilde f_2=0.$$
This shows that ${\mathcal  B}-{\mathcal  A}^H$ is vertical.
\hfill\rule{3mm}{3mm}

\bigskip

Define  $\alpha: P\to {\mathfrak g}\otimes{\mathfrak g }$ and
 $\beta: P\to {\mathfrak g}$
by
$$\alpha(u)=\sum\alpha^{ij}(u)A_i\otimes A_j$$
$$\beta(u)=\sum\beta^k(u) A_k.$$
It is easy to see that $\B^V$ depends only on $\alpha$,  $\beta$ and
 the expression is independent of the choice of basis of
$ {\mathfrak g}$. From the invariance of $\B$ we obtain
\begin{eqnarray*}
\alpha(ug)&=&\left(ad(g)\otimes ad(g)\right) \alpha(u),\\
\beta(ug)&=&ad(g)\beta(u)
\end{eqnarray*}
 for all $u\in P$ and $g\in G$.
\medskip

{\bf C.}
Theorem \ref{theorem-operator-decomposition} has a more directly
 probabilistic version. For this let $\pi:P\to M$ be as before and
for $0\le l<r<\infty$ let $C(l,r;P)$ be the space of continuous paths
$y: [l,r]\to P$ with its usual Borel $\sigma$-algebra. For such
write $l_y=l$ and $r_y=r$. Let $C(*,*;P)$ be the union of such
spaces. It has the standard additive structure under concatenation:
if $y$ and $y'$ are two paths with $r_y=l_{y'}$ and
$y(r_y)=y'(l_{y'})$ let $y+y'$ be the corresponding element in
$C(l_y,r_{y'};P)$.
The {\it basic} $\sigma$-algebra of $C(*,*,P)$ is defined to be the pull
back by $\pi$ of the usual Borel $\sigma$-algebra on $C(*,*;M)$.

Consider the laws $\{\PP_a^{l,r}: 0\le l<r, a\in P\}$ of the
process running from $a$ between times  $l$ and $r$, associated to a
smooth diffusion operator $\B$ on $P$. Assume for simplicity that
the diffusion has no explosion. Thus $\{\PP_{a}^{l,r}, a\in P\}$
is a kernel from $P$ to $C(l,r;P)$. The right action $R_g$ by $g$ in $G$
extends to give a right action, also written $R_g$, of $G$ on $C(*,*,P)$.
Equivariance of $\B$ is equivalent to
$$\PP_{ag}^{l,r}=(R_g)_*\PP_a^{l,r}$$
for all $0\le l\le r$ and $a\in P$. If so $\pi_*(\PP_a^{l,r})$
depends only on $\pi(a)$, $l$, $r$ and gives the law of the
 induced diffusion $\A$ on $M$.
We say that such a diffusion $\B$ is {\it basic} if for all $a\in P$ and
 $0\le l< r<\infty$ the basic $\sigma$-algebra on $C(l,r;P)$ contains
all Borel sets up to $\PP_a^{l,r}$ negligible sets, i.e. for all
$a\in P$ and Borel subsets $B$ of $C(l,r;P)$ there exists a Borel
subset $A$ of $C(l,r,M)$ s.t. $\PP_a(\pi^{-1}(A)\Delta B)=0$.

For paths in $G$ it is more convenient to consider the space
$C_{id}(l,r;G)$ of continuous $\sigma: [l,r]\to G$ with
$\sigma(l)=id$ for $`id'$ the identity element. The corresponding
space $C_{id}(*,*,G)$ has a multiplication
$$C_{id}(s,t;G)\times C_{id}(t,u;G)\longrightarrow C_{id}(s,u;G)$$
$$(g, g')\mapsto g\times g'$$
where $(g\times g')(r)=g(r)$ for $r\in [s,t]$ and
 $(g\times g')(r)=g(t)g'(r)$ for $r\in [t,u]$.

Given probability measures $\Q$, $\Q'$ on $C_{id}(s,t;G)$ and
 $C_{id}(t,u;G)$ respectively this determines a convolution
$\Q*\Q'$ of $\Q$ with $\Q'$ which is a probability measure on
 $C_{id}(s,u;G)$.

\begin{theorem}\label{theorem-kernel-decomposition}
Given the laws  $\{\PP_a^{l,r}: a\in P, 0\le l<r<\infty\}$ of an
equivariant diffusion $\B$ as above with $\A$ strongly cohesive
there exist probability kernels $\{\PP_a^{H,l,r}: a\in P\}$ from
$P$ to $C(l,r;P)$, $0\le l < r <\infty$ and $\Q_y^{l,r}$,
 defined $\PP^{l,r}$  a.s. from $C(l,r,P)$ to $C_{id}(l,r;G)$ such
that
\begin{enumerate}
\item[(i)]
$\{\PP_a^{H,l,r}: a\in P\}$ is equivariant, basic and determining a 
strongly cohesive generator.
\item[(ii)]
$y\mapsto \Q_y^{l,r}$ satisfies
$$\Q_{y+y'}^{l_y, r_{y'}}=\Q_y^{l_y, r_y}*\Q_{y'}^{l_{y'},r_{y'}}$$
for $\PP^{l_{y}, r_y}\otimes \PP^{l_{y'},r_{y'}}$ almost all
$y$, $y'$ with $r_y=l_{y'}$.
\item[(iii)]
For $U$ a Borel subset of $C(l,r,P)$,
$$\PP_a^{l,r}(U)
=\int\int\chi_{U}(y_\cdot\cdot  g_\cdot) \Q_y^{l,r}(dg)
\PP_a^{H,l,r}(dy).$$
\end{enumerate}
The kernels $\PP_a^{H,l,r}$ are uniquely determined as are the
$\{\Q_y^{l,r}: y\in \R\}$, $\PP_a^{H,l,r}$ a.s. in $y$ for all $a$
in $P$. Furthermore $\Q_y^{l,r}$ depends on $y$ only through its
projection $\pi(y)$ and its initial point $y_l$.
\end{theorem}

{\it Proof.}
Fix $a$ in $P$ and let $\{b_t: l\le r\le t\}$ be a process with law
$\PP_a^{l,r}$.
By Theorem \ref{theorem-operator-decomposition} we can assume that
 $b_\cdot$ is given by an  s.d.e. of the form
\begin{equation}
\label{bt}
db_t=\tilde X(b_t)\circ dB_t+\tilde X^0(b_t)dt
+A(b_t)\circ d\beta_t+V(b_t)dt
\end{equation}
where $\tilde X: P\times \R^p\to TP$ is the horizontal lift of
some $X: M\times \R^p\to E$, $\tilde X^0$ is the horizontal lift
of a vector field $X^0$ on $M$, while $A: P\times \R^1\to TP$ and
the vector field $V$ are vertical and determine $\B^V$. Here
$B_\cdot$ and $\beta_\cdot$ are independent Brownian motions on
$\R^p$ and $\R^q$
 respectively, some $q$, and we are using the semi-connection on $P$ induced
by $\B$ as in Theorem \ref{theorem:connection}.

 Let $\{\tilde x_t: l\le t \le r\}$ satisfy
\begin{equation}
\label{sde-on-p}
\begin{aligned}
d\tilde x_t&=\tilde X(\tilde x_t)\circ dB_t
+\tilde X^0(\tilde x_t)dt\\
\tilde x_l&=a
\end{aligned}
\end{equation}
so $\tilde x_\cdot$ is the horizontal lift of
 $\{\pi(b_t): l\le t\le r \}$. Then there is a unique continuous
 process $\{g_t: l\le t\le r\}$ in $G$ with $g_l=id$ such that
$$\tilde x_tg_t=b_t.$$

We have to analyse $\{g_t: l\le t \le r\}$. Using local
trivialisations of $\pi:P\to M$ we see it is a semi-martingale. As in
 \cite{Kobayashi-Nomizu-I}, Proposition 3.1 on page 69,
$$db_t=TR_{g_t}(\circ d\tilde x_t)+A^{g_t^{-1}\circ dg_t}(b_t)$$
giving $$\tilde \omega(\circ db_t)=\tilde \omega
\left(A^{g_t^{-1}\circ dg_t}(b_t)\right)=g_t^{-1}\circ dg_t$$
for any smooth connection form $\tilde \omega: P\to {\mathfrak g}$
on $P$ which vanishes on $H^ETP$. Thus
\begin{equation}
\begin{aligned}
dg_t&=TL_{g_t} \tilde \omega\big(A(\tilde x_tg_t)\circ
d\beta_t+ V(\tilde x_tg_t)dt\big)\\
g_l&=id, \qquad l\le t \le r.
\end{aligned}
\end{equation}
For $y\in C(l,r:P)$ let $\{g_t^y: l\le t \le r\}$ be the solution of
\begin{equation}
\begin{aligned}
dg_t^y&=TL_{g_t^y} \tilde \omega\big(A(y_tg_t^y)\circ d\beta_t+
 V(y_tg_t^y)dt\big)\\
g_l^y&=id
\end{aligned}
\end{equation}
(where the Stratonovich equation is interpreted with `$dy_td\beta_t=0$').
Since $\beta_\cdot$ and $B_\cdot$ and hence $\beta_\cdot$ and
 $\tilde x_\cdot$ are independent we see $g=g^{\tilde x}$ almost surely.
 For a discussion of some technicalities concerning skew products,
 see \cite{Taylor92}.

For $y_\cdot$ in $C(*,*;P)$ let $\{h(y)_t: l_y\le t\le r_y\}$ be
the horizontal lift of $\pi(y)_\cdot$, starting at $y_{l_y}$. This
exists for almost all $y$ as can be seen either by the extension
of It\^o's result to general principal bundles, e.g. using
(\ref{sde-on-p}), or by the existence of measurable sections using
the fact that $\A^H$ is basic. Define $\PP_a^{H,l,r}$ to be the
law of $\tilde x_\cdot $ above and $Q_y^{l,r}$ to be that of
$g^{h(y)}$. Clearly conditions (i) is satisfied.

To check (ii) take $y$ and $y'$ with $r_y=l_{y'}$. Then
$$h(y+y')=h(y)+h(y') \left(g_{r_y}^{h(y)}\right)^{-1},$$
writing $y=h(y)g^{h(y)}$ and  $y'=h(y')g^{h(y')}$.
For $r_y\le t\le r_{y'}$ this shows
$$(y+y')_t=h(y')_t\left(g_{r_y}^{h(y)}\right)^{-1} g_t^{h(y+y')}.$$
But $(y+y')_t=y_t'=h(y')_tg_t^{h(y')}$ and so we have
$g_t^{h(y+y')}=g_{r_y}^{h(y)}g_t^{h(y')}$ for $t\ge r_y$, giving
$g^{h(y+y')}=g^{h(y)}\times g^{h(y')}$ almost surely. This proves
(ii).

For uniqueness suppose we have another set of probability measures
denoted $\tilde \Q_y^{l,r}$ and $\tilde P_a^{H,l,r}$ which satisfy
(i), (ii), (iii). Since $\{\tilde \PP_a^{H,l,r}\}_a$ is
equivariant and induces $\A$ on $M$ we can apply the preceding
argument to it in place of $\{\PP_a^{l,r}\}_a$. However since it
is basic the term involving $\beta$ in the stochastic differential
equation (\ref{sde-on-p}) must vanish. Since it is also strongly
cohesive the vertical part $V$ must vanish also and we have
$\tilde \PP_a^{H,l,r} =\PP_a^{H,l,r}$. On the other hand in the
decomposition
 $b_t=\tilde x_t g_t^{\tilde x_t}$ the law of $g_\cdot^{\tilde x}$ is
 determined  by those of $b_\cdot$ and $\tilde x_\cdot$ but $\Q_y^{l,r}$
 is the conditional law of $g_\cdot^{\tilde x_\cdot}$ given $\tilde x=y_\cdot$
and so is uniquely determined  as described.
\hfill\rule{3mm}{3mm}
\medskip

In fact $\Q_{y}^{l,r}$ is associated to the time dependent generator
which at $g\in G$ and $t\in[l,r]$ is
 $\sum \alpha^{ij}(h(y)_tg) {\mathcal  L}_{A_i}{\mathcal  L}_{A_j}
  +\sum\beta^k(h(y)_tg){\mathcal  L}_{A_k}$ for $\alpha^{ij}$ and
$\beta^k$ as defined in Theorem \ref{theorem-operator-decomposition}
while $\PP^{H,l,r}$ is associated to $\A^H$.

\section{Stochastic flows and derivative flows}
\label{section-stochastic-flows}

{\bf A. Derivative flows.}
 Let ${\mathcal A}$ on $M$ be given in H\"ormander form
$${\mathcal  A}
 ={1\over 2}\sum_{j=1}^m {\mathcal  L}_{X^j}{\mathcal  L}_{X^j}
+{\mathcal  L}_{A}$$ for some vector fields
$\displaystyle{X^1,\dots X^m}$, $A$. As before let
$E_x=span\{X^1(x),\dots, X^m(x)\}$ and assume
 $dim E_x$ is constant, $p$,
say, giving a sub-bundle $E\subset TM$. The
$\displaystyle{X^1(x),\dots, X^m(x)}$ determine a vector bundle map
of the trivial bundle $\underline {\R}^m$
$$X:\underline {\R}^m\longrightarrow   TM$$
with $\sigma^{\mathcal  A}=X(x)X(x)^*.$ We can, and will, consider
$X$ as a map $\displaystyle{X: \underline{\R}^m\to E}$.

As such it determines (a) a Riemannian metric
 $\{\langle, \rangle_x: x \in M\}$ on $E$ (the same as that determined
by $\sigma^A$) and (b) a metric connection $\breve \nabla$ on $E$
uniquely defined by the requirement that for each $x$ in $M$,
$$\breve \nabla_vX(e)=0$$
for all $v\in T_xM$ whenever $e$ is orthogonal to the kernel of $T_xM$.
 Then for any differentiable section $U$ of $E$,
\begin{equation}
\label{def-semi-connection} \breve \nabla_vU=Y(x)
d\left(Y_\cdot(U(\cdot))\right)(v), \qquad v \in T_xM,
\end{equation}
 where $Y$ is the $\R^m$ valued 1-form on $M$ given by
$$\langle Y_x(v),e\rangle_{\R^m}=\langle X(x)(e), v \rangle_x,
\qquad e\in \R^m, v\in E_x, x\in M$$ e.g.
\cite{Elworthy-LeJan-Li-book} where it is referred to as the
LeJan-Watanabe connection in this context. By a theorem of
 Narasimhan and Ramanan \cite{Narasimhan-Ramanan61} all metric connections
 on $E$ arise this way, see \cite{Quillen}, \cite{Elworthy-LeJan-Li-book}.

 For $\displaystyle{\{B_t: 0\le t<\infty\}}$ a Brownian motion on $\R^m$,
 the  stochastic differential equation
\begin{equation}\label{sde}
dx_t= X(x_t)\circ dB_t +A(x_t)dt
\end{equation}
determines a Markov process with differential generator ${\mathcal
A}$. Over each solution $\displaystyle{\{x_t: 0\le t<\rho\}}$,
where $\rho$ is the explosion time,  there is a
 `derivative' process $\displaystyle{\{v_t: 0\le t<\rho\}}$ in $TM$ which
we can write as $\displaystyle{\{T\xi_t(v_0):0\le t<\rho\}}$ with
$\displaystyle{T\xi_t: T_{x_0}M\to T_{x_t}M}$ linear. This would
be the
 derivative  of the flow $\displaystyle{\{\xi_t: 0\le t<\rho\}}$
of the stochastic differential equation when the stochastic
differential
 equation is strongly complete. In general it is given by a
stochastic differential equation  on the tangent bundle $TM$, or
equivalently by a covariant equation along $\displaystyle{\{x_t:
0\le t<\rho\}}$:
$$Dv_t=\nabla X(v_t)\circ dB_t+\nabla A(v_t)dt$$
with respect to any torsion free connection. Take $P$ to be the
linear frame
 bundle $GL(M)$ of $M$, treating $\displaystyle{u\in GL(M)}$ as an
 isomorphism $\displaystyle{u:\R^n\to T_{\pi(u)}M}$. For $u_0\in GLM$ we
obtain a process $\displaystyle{\{u_t: 0\le t<\rho\}}$ on $GLM$ by
$$u_t=T\xi_t\circ u_0.$$ Let ${\mathcal B}$ be its differential
generator. Clearly it is
 equivariant and a lift of ${\mathcal  A}$.
\medskip

A proof of the following in the context of stochastic flows, is
given later. For $w\in E_x$, set
\begin{equation}
\label{Zw} Z^w(y)=X(y)Y(x)(w).
\end{equation}

\begin{theorem}\label{Derivative-connection}
The semi-connection $\nabla$ induced by ${\mathcal  B}$ is the
adjoint connection of the LeJan-Watanabe connection $\breve\nabla$
determined by $X$,
 as defined by  (\ref{def-semi-connection}), \cite{Elworthy-LeJan-Li-book}.
 Consequently $\nabla_wV=L_{Z^w}V$ for any vector field $V$ and $w\in E$ also
$\nabla_{V(x)} Z^w$ vanishes  if $w\in E_x$.
\end{theorem}

In the case of the derivative flow the $\alpha$, $\beta$ of Theorem
  \ref{theorem-operator-decomposition} have an explicit expression:
for $u\in GLM$,
\begin{equation}
\label{derivative-symbols}
 \left\{
\begin{aligned}
\alpha(u)&={1\over 2}\sum \left(u^{-1}(-)\breve \nabla_{u(-)}X^p\right)
\otimes \left(u^{-1}(-)\breve\nabla_{u(-)} X^p\right)\\
\beta(u)&=-{1\over 2}\sum u^{-1} \breve \nabla_{\breve \nabla_{u(-)}X^p}X^p
- {1\over 2}u^{-1} \hbox{Ric }^{\#}{u(-)}.
\end{aligned}\right.
\end{equation}
Here $\breve R$ is the curvature tensor for $\breve \nabla$ and
$\breve Ric^{\#}: TM\to E$  the Ricci curvature defined by
$\breve Ric^{\#}(v)=\sum_{j=1}^p \breve R(v, e^j)e^j$, $v\in
T_xM$.

\medskip

Equivariant operators on $GLM$ determine operators on associated
bundles, such as $\wedge^q TM$. If the original operator was
vertical this turns out to be a zero order operator (as is shown
in \cite{Elworthy-LeJan-Li-principal} for general principal bundles) and
in the case of $\wedge^q TM$ these operators are the generalized
Weitzenbock curvature operators described in
\cite{Elworthy-LeJan-Li-book}. In particular for differential
1-forms the operator is $\phi\mapsto \phi(Ric^{\#}-)$. To see
this, as an illustrative example, given a 1-form $\phi$ define
$\tilde \phi: GLM\to L(\R^n;\R)$ by $\tilde \phi(u)=\phi_{\pi u}u$
so $\tilde\phi(ug)=\phi_{\pi u}(ug-)$. Then
\begin{eqnarray*}
L_{A_j^*}(\tilde\phi)(u)&=& {d\over dt} \tilde \phi(u\cdot
e^{A_jt})
\vert_{t=0}\\
&=& {d\over dt}  \phi_{\pi u}(u\cdot e^{A_jt})
\vert_{t=0}\\
&=& \phi_{\pi u}(u A_j-)=\tilde \phi(u) (A_j-).
\end{eqnarray*}
Iterating we have
\begin{eqnarray*}
\B^V(\tilde\phi)(u)&=& \sum_{i,j}\alpha^{i,j}(u)
  \phi_{\pi u}(u A_jA_i-)+\sum_{k}\beta^{k}(u)
  \phi_{\pi u}(u A_k-)\\
&=&-{1\over 2}\tilde\phi(u)(u^{-1}Ric^{\#}(u-))
\end{eqnarray*}
as required, by using the map $gl(n)\otimes gl(n)\to gl(n)$,
 $ S\otimes T\mapsto S\circ T$, and equation
 (\ref{derivative-symbols}).

\medskip

{\bf B. Stochastic flows.}
In fact Theorem~\ref{Derivative-connection} can be understood in the more
general context of stochastic flows as diffusions on the diffeomorphism
groups. For this assume that $M$ is compact and for $r\in\{1,2,\dots\}$ and
  $s>r+\hbox{dim}(M)/2$ let
$\displaystyle{{\mathcal D}^s={\mathcal D}^sM}$ be the
 $C^\infty$ manifold of  diffeomorphisms of $M$ of  Sobolev class $H^s$,
(for example see  Ebin-Marsden \cite{Ebin-Marsden} or
  Elworthy \cite{Elworthy-book}.) Alternatively we could take the space
$\D^\infty$ of $C^\infty$ diffeomorphisms with differentiable
structure as in \cite{Kriegl-Michor}. Fix a base point $x_0$ in
$M$ and let $\pi: \D^s\to M$ be  evaluation at $x_0$. This makes
$\D^s$ into a principal bundle over $M$ with group the manifold
$\D_{x_0}^s$ of $H^s$- diffeomorphisms $\theta$ with
$\theta(x_0)=x_0$, acting on the right by composition (although
the action of $\D^{s+r}$ is only $C^r$, for $r=0,1,2,\dots$).

Let  $\{\xi_t^s: 0\le s\le t<\infty\}$ be the flow of (\ref{sde})
starting at time $s$. Write $\xi_t$ for $\xi^0_t$.
The more general case allowing for infinite
dimensional noise is given in \cite{Elworthy-LeJan-Li-principal}.
We define probability measures $\{\PP_{\theta}^{s,t}: \theta\in
\D^s\}$ on $C([s.t]; M)$ be letting $\PP_\theta^{s,t}$ be the law
of $\{\xi_r^s\circ \theta: s\le r\le t\}$ (These correspond to the
diffusion process on $\D^s$ associated to the right-invariant
stochastic differential equation on $\D^s$ satisfied by $\{\xi_t:
0\le t<\infty\}$ as in \cite{Elworthy-book}.) These are
equivariant and project by $\pi$ to the laws given by the
stochastic differential equation on $M$. Assuming that these give a
strongly cohesive diffusion on $M$ we are essentially in the
situation of Theorem \ref{theorem-kernel-decomposition}.

Let $K(x): \R^m\to \R^m$ be the orthogonal projection onto the kernel
of $X(x)$, each $x\in M$.  set $K^\perp(x)=id-K(x)$. Consider the
$\D^\infty$-valued process $\{\theta_t: 0\le t<\infty\}$ given
by (or as the flow of)
\begin{equation}
\label{diffeo-1}
d\theta_t(x)=X(\theta_t(x))K^\perp(\theta_t(x_0))\circ dB_t
   +X(\theta_t(x))Y(\theta_t(x_0))A(\theta_t(x_0))
\end{equation}
for given $\theta_0$ in $\D^\infty$ and, define a
 $\D_{x_0}^\infty$-valued process $\{g_t: 0\le t<\infty\}$ by
\begin{eqnarray}\label{process-group}
dg_t&=&T\theta_t^{-1}\left\{X(\theta_tg_t-)K(\theta_tx_0)\circ dB_t\right.\\
&&
\nonumber\left.
+A(\theta_tg_t-)dt-X(\theta_tg_t-)Y(\theta_tx_0)A(\theta_tx_0)dt\right\}\\
g_0&=&id.
\nonumber
\end{eqnarray}
Set $x_t^\theta=\xi_t(\theta_0(x_0))$. Note that
 $\pi(\theta_t)=\theta_t(x_0)=x_t^\theta$ since
$$X(\theta_t(x_0))K^\perp(\theta_t(x_0))=X(\theta_t(x_0))$$
and
$$X(\theta_t(x_0))Y(\theta_t(x_0))A(\theta_t(x_0))=A(\theta_t(x_0)).$$
Thus $\{\theta_t: 0\le t<\infty\}$ is a lift of
 $\{x_t^\theta, 0\le t<\infty\}$. It can be considered to be driven by
 the `relevant noise',
(from the point of view of $\xi_\cdot(\theta_0(x_0))$, i.e. by the
Brownian motion $\tilde B_\cdot$ given by
$$\tilde B_t=\int_0^t\tilde \paral(x_\cdot^\theta)_s^{-1}
K^\perp(x_s^\theta)dB_s$$
where $\{\tilde{\paral}(x_\cdot^\theta), 0\le s<\infty\}$
is parallel translation along $
x_\cdot^\theta$ with respect to the connection on the trivial bundle
$M\times \R^m\to M$ determined by $K$ and $K^\perp$, so that
$$\tilde \paral(x_\cdot^\theta)_s: \R^m\to \R^m$$
is orthogonal and maps the kernel of $X(\theta_\cdot(x_0))$
 onto the kernel of $X(x_s^\theta)$ for $0\le s<\infty$,
 see \cite{Elworthy-LeJan-Li-book}(chapter 3).

Correspondingly there is the `redundant noise', the Brownian motion
$\{\beta_t:~0\le t<\infty\}$ given by
$$\beta_t=\int_0^t \tilde
 \paral(x_\cdot^\theta)_s^{-1} K(x_s^\theta)dB_s.$$
Then, as shown in  \cite{Elworthy-LeJan-Li-book}(chapter 3),
\begin{enumerate}
\item[(i)]
$\tilde B_\cdot$ has the same filtration as $\{x_s^\theta: 0\le s<\infty\}$
\item[(ii)]
$\beta_\cdot$ and $\tilde B_\cdot$ are independent
\item[(iii)]
$dB_t=\tilde{\parals_t} d\beta_t +\tilde{\parals_t} d\tilde B_t$.
\end{enumerate}

We wish to see how $g_\cdot$ is driven by $\beta_\cdot$. For this observe
$$\int_0^t K(x_s^\theta)\circ dB_s=
\int_0^t K(x_s^\theta) dB_s +\int_0^t \Lambda(x_s^\theta)ds$$
for $\Lambda: M \to \R$ given by the Stratonovich correction term.
By (iii) $$\int_0^t K(x_s^\theta)dB_s=\int_0^t \tilde{\parals_s} d\beta_s
=\int_0^t \tilde{\parals_s} \circ d\beta_s$$
 since $\tilde{\parals_\cdot}$
is independent of $\beta$ by (i) and (ii). Thus
equation~(\ref{process-group}) for $g_\cdot$ can be written as
\begin{eqnarray*}
dg_t&=&T\theta_t^{-1}\left\{X(\theta_tg_t-) \tilde{\paral}
(\theta_\cdot(x_0))_t
 \circ d\beta_t +  X(\theta_tg_t-) \Lambda(\theta_t(x_0))dt \right.\\
&&\left.
+A(\theta_tg_t-)dt-X(\theta_tg_t-)Y(\theta_tx_0)A(\theta_tx_0)dt\right\}
\end{eqnarray*}
and if we define
\begin{eqnarray*}
dg_t^y&=&Ty_t^{-1}\left\{X(y_tg_t-) \tilde{\paral}(y_\cdot(x_0))_t
 \circ d\beta_t + X(y_tg_t-) \Lambda(y_t(x_0))dt \right.\\
&&\left.
+A(y_tg_t-)dt-X(y_tg_t-)Y(y_tx_0)A(y_tx_0)dt\right\}\\
g_0&=&id
\end{eqnarray*}
for any continuous $y: [0,\infty)\to \D^\infty$, we see, by the
independence of $\beta$ and $\theta$ that
$g_\cdot=g_\cdot^\theta$.

By It\^o's formula  on $\D^s$, for $x\in M$,
$$d(\theta_t g_t(x_0))=(\circ d\theta_t) (g_t(x))+
T\theta_t (\circ dg_t^\theta(x)).$$
Now
\begin{eqnarray*}
T\theta_t (\circ dg_t^\theta(x))
&=&\left\{X(\theta_tg_t(x))K(\theta_tx_0)\circ dB_t\right.\\
&&\left.
+A(\theta_tg_t(x))dt-X(\theta_tg_t(x))Y(\theta_tx_0)A(\theta_tx_0)dt\right\}
\end{eqnarray*}
and so by (\ref{diffeo-1}) we see that $\theta_tg_t=\xi_t\circ \theta_0$,
 a.s.

Taking $\theta_0=id$ we have

\begin{proposition}
The flow $\xi_\cdot$ has the decomposition
$$\xi_t=\theta_t g_t^{\theta_\cdot}, \qquad
0\le t<\infty$$ for $\theta$ and $g_\cdot^\theta\equiv g_\cdot$
given by (\ref{diffeo-1}) and (\ref{process-group}) above. For
almost all $\sigma: [0,\infty)\to M$ with $\sigma(0)=x_0$ and
bounded measurable $F: C(0,\infty; \D^\infty)\to \R$
$$\E\left\{ F(\xi_\cdot) \left\vert
\xi_\cdot(x_0)=\sigma\right.\right\}
=\E\left\{ F(\tilde \sigma g_\cdot^{\tilde \sigma}) \right\}$$
where $\tilde \sigma: [0,\infty)\to \D^\infty$ is the horizontal lift of
$\sigma$ with $\tilde \sigma(0)=id$.
\end{proposition}
To define the `horizontal lift' above we can use the fact, from (i)
above, that $\theta_\cdot$ has the same filtration as $\xi_\cdot(x_0)$
and so furnishes a lifting map.

In terms of the semi-connection induced on $\pi: \D^s\to M$ over
$E$, from above, by uniqueness or directly, we see the horizontal
lift
\begin{eqnarray*}
h_\theta&:& E_{\theta(x_0)}\longrightarrow T_\theta\D^s\\
h_\theta(v)&:&M\longrightarrow TM
\end{eqnarray*}
is given by $h_\theta(v)=X(\theta(x))Y(\theta(x_0))v$ and the
horizontal lift $\tilde \sigma_\cdot$ from $\tilde \sigma_{0}$ of
a $C^1$ curve $\sigma$ on $M$ with $\tilde \sigma_0(x_0)=\sigma_0$
and $\dot \sigma(t)\in E_{\sigma(t)}$, all $t$,  is given by
$${d\over dt}\tilde \sigma_t=X(\tilde\sigma_t-)Y(\sigma_t)\dot \sigma_t$$
for $\tilde \sigma_0=id$. The lift is the solution flow of the
differential equation
$$\dot y_t=Z^{\dot \sigma}(y_t)$$
on $M$.
\medskip

For each frame $u:\R^n\to T_{x_0}M$ there is a homomorphism of principal
bundles
\begin{equation}
\label{principal-bundle-homo}
\begin{array}{ll}
\D^s&\to GL(M)\\
\theta &\mapsto T_{x_0}\theta\circ u.
\end{array}
\end{equation}
This sends $\{\xi_t: t\ge 0\}$ to the derivative process
$T_x\xi_t\circ u$. (If the latter satisfies the strongly cohesive
condition we could apply our analysis to this submersion $\D^s\to
GLM$ and get another decomposition of $\xi_\cdot$.)

\medskip

 Results in Kobayashi-Nomizu \cite{Kobayashi-Nomizu-I}
(Proposition 6.1 on page 79) apply to the homomorphism $\D^s\to
GL(M)$ of (\ref{principal-bundle-homo}). This gives a relationship
between the curvature and holonomy groups of the semi-connection $\hat
\nabla$ on $GLM$ determined by the derivative flow  and those of
the connection induced by the diffusion on
$\D^s\stackrel{\pi}{\longrightarrow} M$. It also shows  that the
horizontal lift $\{\tilde x_t: t\ge 0\}$ through $u$ of $\{x_t:
t\ge 0\}$ to $GL(M)$ is just $T_{x_0}\theta_t\circ u$ for
$\{\theta_t: t\ge 0\}$ the flow given by (\ref{diffeo-1})
 with $\theta_0=id$, i.e. the solution flow of the stochastic
 differential equation
$$dy_t=Z^{\circ dx_t}(y_t).$$
From this and Lemma 1.3.4 of \cite{Elworthy-LeJan-Li-book} we see
that $\hat \nabla$ is the adjoint of the LeJan-Watanabe connection
determined by the flow, so proving Theorem
\ref{Derivative-connection} above. However the present
construction applies with GLM replaced by any natural bundle over
$M$ (e.g. jet bundles, see Kolar-Michor-Slovak
\cite{Kolar-Michor-Slovak}), to give semi-connections on these
bundles.

\bigskip

{\tiny \hskip-\parindent K. D. ELWORTHY, \hskip4pt MATHEMATICS
INSTITUTE,  WARWICK UNIVERSITY, COVENTRY CV4 7AL, UK}
\newline
{\tiny \hskip-\parindent Y. LE JAN, \hskip4pt D\'EPARTMENT DE
MATH\'EMATIQUE, UNIVERSIT\'E PARIS SUD, 91405 ORSAY,
FRANCE}\newline {\tiny \hskip-\parindent XUE-MEI LI, \hskip4pt
DEPARTMENT OF COMPUTING AND MATHEMATICS, THE NOTTINGHAM TRENT
UNIVERSITY, NOTTINGHAM NG7 1AS, UK.} {\scriptsize e-mail address:
xuemei.li@ntu.ac.uk}

\end{document}